\input amstex 
\documentstyle{amsppt} 
\magnification=1200
\pagewidth{6.5truein}
\pageheight{8.7truein}
\nologo
\define\Hom{\operatorname{Hom}}
\define\colim{\operatornamewithlimits{colim}}
\input xy
\xyoption{all}
\NoRunningHeads

\topmatter

\title Mapping cylinders and the Oka Principle \endtitle
\author Finnur L\acuteaccent arusson \endauthor
\affil University of Western Ontario \endaffil
\address Department of Mathematics, University of Western Ontario,
         London, Ontario N6A~5B7, Canada \endaddress
\email   larusson\@uwo.ca \endemail
\thanks The author was supported in part by the Natural Sciences and 
Engineering Research Council of Canada. \endthanks
\thanks 16 September 2004. \endthanks

\abstract
We apply concepts and tools from abstract homotopy theory to complex
analysis and geometry, continuing our development of the idea that the
Oka Principle is about fibrancy in suitable model structures.  We
explicitly factor a holomorphic map between Stein manifolds through
mapping cylinders in three different model structures and use these
factorizations to prove implications between ostensibly
different Oka properties of complex manifolds and holomorphic maps.
We show that for Stein manifolds, several Oka properties coincide and
are characterized by the geometric condition of ellipticity.  Going
beyond the Stein case to a study of cofibrant models of arbitrary
complex manifolds, using the Jouanolou Trick, we obtain a geometric
characterization of an Oka property for a large class of manifolds,
extending our result for Stein manifolds.  Finally, we prove a converse
Oka Principle saying that certain notions of cofibrancy for manifolds
are equivalent to being Stein.

\endabstract

\endtopmatter

\document

\subhead Introduction \endsubhead
In this paper, we apply concepts and tools from abstract homotopy theory
to complex analysis and geometry, based on the foundational work in
\cite{L2}, continuing our development of the idea that the Oka Principle
is about fibrancy in suitable model structures.  A mapping cylinder
in a model category is an object through which a given map can be
factored as a cofibration followed by an acyclic fibration (or sometimes
merely an acyclic map).  We explicitly factor a holomorphic map between
Stein manifolds through mapping cylinders in three different model
structures.  We apply these factorizations to the Oka Principle, mainly
to prove implications between ostensibly different Oka properties of
complex manifolds and holomorphic maps.  We show that for Stein manifolds,
several Oka properties coincide and are characterized by the geometric
condition of ellipticity.  We then move beyond the Stein case to a study
of cofibrant models of arbitrary complex manifolds (this involves the
same sort of factorization through a mapping cylinder as before).  Using
the so-called Jouanolou Trick, we obtain a geometric characterization of
an Oka property for a large class of manifolds, extending our result
for Stein manifolds.  Finally, we prove a \lq\lq converse Oka
Principle\rq\rq\ saying that certain notions of cofibrancy for manifolds
are equivalent to being Stein.  This confirms and makes precise the
impression that the Stein property is dual to the Oka properties that
have been proved equivalent to fibrancy.

\subhead Factorization in the intermediate structure \endsubhead
We shall work in the simplicial category $\frak S$ of prestacks on the
simplicial site $\Cal S$ of Stein manifolds.  The category of complex
manifolds is simplicially embedded in $\frak S$.  There are at least six
interesting simplicial model structures on $\frak S$: the projective,
intermediate, and injective structures with respect to either the
trivial topology on $\Cal S$ or the topology in which a cover, roughly
speaking, is a cover by Stein open subsets.  These structures may be
used to study the Oka Principle for complex manifolds and holomorphic
maps.  We sometimes refer to notions associated with the former topology
as {\it coarse} and the latter as {\it fine}.  We often omit the word
{\it fine}, as the fine structures are more important than the coarse
ones.  For more background, see \cite{L2}, where the intermediate
structures were first introduced, and the references there.

One of the axioms for a model category states that any map can be
factored (by no means uniquely) as a cofibration followed by an acyclic
fibration.  The general verification of this for a category of prestacks
is quite abstract and not always useful in practice.  We start by
producing an explicit factorization in the intermediate structures ---
they have the same cofibrations and acyclic fibrations ---  for a
holomorphic map $f:R\to S$ between Stein manifolds.  We assume that $R$
is finite-dimensional, meaning that the (finite) dimensions of its
connected components are bounded (by convention, we take manifolds to
be second countable but not necessarily connected).

Let $\phi:R\to\Bbb C^r$ be a holomorphic embedding of $R$ into
Euclidean space.  Let $M$ be the Stein manifold $S\times\Bbb C^r$.
Then $f$ can be factored as $R @>j>> M @>\pi>> S$, where $j=(f,\phi)$
is an embedding and $\pi$ is the projection $(x,z) \mapsto x$.  In the
language of \cite{L2}, $j$ is a Stein inclusion, so it is an
intermediate cofibration.  Since $\Bbb C^r$ is elliptic, $\pi$ is an
intermediate fibration: it is the pullback of the constant map from
$\Bbb C^r$ by the constant map from $S$.  Also, $\pi$ is pointwise
acyclic since $\Bbb C^r$ is holomorphically contractible.  Following
standard terminology in topology, we refer to $M$ as a {\it mapping
cylinder} for $f$.

\subhead Interpolation implies approximation \endsubhead
Complex manifolds and, as introduced in \cite{L2}, holomorphic maps can
possess various Oka properties.  First, there are the parametric Oka
property (POP) and the strictly stronger parametric Oka property with
interpolation on a submanifold (POPI).  For manifolds and those
holomorphic maps that are topological fibrations (in the sense of Serre
or Hurewicz: these are equivalent for continuous maps between smooth
manifolds), POP, called the weak Oka property in \cite{L2}, is equivalent
to projective fibrancy, and POPI, called {\it the} Oka property in
\cite{L2}, is equivalent to intermediate fibrancy.

There is also the parametric Oka property with approximation on
a holomorphically convex compact subset (POPA): this property has not been
explicitly generalized from manifolds to maps until now; the natural notion
appears below.  Each parametric property has a {\it basic} version
(BOP, BOPI, BOPA), in which the parameter space is simply a point.  A
disc satisfies POP but neither BOPI nor BOPA.  For more background, see
\cite{L2} and the survey \cite{F2}.

Approximation on holomorphically convex compact subsets is important
in most areas of complex analysis, including basic proofs in Oka
theory.  It is therefore of interest that POPI implies POPA.
The following argument is short and easy and makes use of the mapping
cylinder introduced above.  It seems reasonable to view POPA as a
reflection of a special case of POPI, but not as a possibly distinct
notion of intermediate fibrancy in its own right.  We do not
know whether POPA is strictly weaker than POPI: examples in this area
are hard to come by.

Let $k:X\to Y$ be a holomorphic map satisfying POPI.  Let $P$ be a
finite polyhedron with a subpolyhedron $Q$.  Let $S$ be a Stein
manifold and $K$ be a holomorphically convex compact subset of $S$ with
a neighbourhood $U$.  Let $h:S\times P\to X$ be a continuous map such that
\roster
\item"{(a)}"  $h(\cdot,p)$ is holomorphic on $U$ for all $p\in P$,
\item"{(b)}" $h(\cdot,q)$ is holomorphic on $S$ for all $q\in Q$, and
\item"{(c)}" $k\circ h(\cdot,p)$ is holomorphic on $S$ for all $p\in P$.
\endroster
Let $d$ be a metric on $X$ compatible with the topology of $X$ and let
$\epsilon>0$.  To show that $k$ satisfies POPA we need to find a
continuous deformation $H:S\times P\times I\to X$, where $I$ is the unit
interval $[0,1]$, such that
\roster
\item $H(\cdot,\cdot,0)=h$,
\item $H(\cdot,p,1)$ is holomorphic on $S$ for all $p\in P$,
\item $H(\cdot,q,t)=h(\cdot,q)$ for all $q\in Q$ and $t\in I$,
\item $k\circ H(\cdot,p,t)=k\circ h(\cdot,p)$ all $p\in P$ and
$t\in I$, and
\item $d(H(x,p,t),h(x,p))<\epsilon$ for all $(x,p,t)\in K\times P\times I$.
\endroster

Let $R\subset U$ be a finite-dimensional Stein neighbourhood of $K$.
The inclusion $f:R\to S$ factors through the mapping cylinder
$M=S\times\Bbb C^r$ introduced above as $R @>j>> M @>\pi>> S$.  Define
$g=h(\pi(\cdot),\cdot):M\times P\to X$. Then $g$ is continuous;
$g(\cdot,p)$ is holomorphic on $R$ (more precisely, $g(j(\cdot),p)$ is
holomorphic) for all $p\in P$; $g(\cdot,q)$ is holomorphic on $M$ for
all $q\in Q$; and $k\circ g(\cdot,p)$ is holomorphic on $M$ for all
$p\in P$.  Since $k$ satisfies POPI, there is a continuous deformation
$G:M \times P\times I\to X$ such that $G(\cdot,\cdot,0)=g$;
$G(\cdot,p,1)$ is holomorphic on $M$ for all $p\in P$;
$G(\cdot,q,t)=g(\cdot,q)$ for all $q\in Q$ and $t\in I$; and finally
$G(j(\cdot),p,t)=g(j(\cdot),p)$ and $k\circ G(\cdot,p,t)=k\circ
g(\cdot,p)$ for all $p\in P$ and $t\in I$.

Since $K$ is holomorphically convex in $S$, the embedding $\phi:R\to\Bbb C^r$
can be uniformly approximated on $K$ as closely as we wish by a holomorphic
map $\psi:S\to\Bbb C^r$.  Precomposing the deformation $G$ in its first
argument by the holomorphic section $\sigma=(\text{id}_S,\psi):S\to M$ of
$\pi$, we obtain a continuous map $H:S\times P\times I\to X$ such that
$$H(\cdot,\cdot,0)=G(\sigma(\cdot),\cdot,0)=g(\sigma(\cdot),\cdot)=
h(\pi\circ\sigma(\cdot),\cdot)=h,$$
so \therosteritem1 holds.  Verifying \therosteritem2, \therosteritem3, and
\therosteritem4 is straightforward.  Finally, for $(x,p,t)$ in the compact
space $K\times P\times I$,
$$\aligned d(H(x,p,t), h(x,p)) &= d(G(\sigma(x),p,t), g(j(x),p)) \\
&=d(G((x,\psi(x)),p,t), G((x,\phi(x)),p,t)) \endaligned$$
is as small as we like.  Let us record what we have just proved.

\proclaim{Theorem 1}  For holomorphic maps between complex manifolds,
the parametric Oka property with interpolation implies the parametric
Oka property with approximation.
\endproclaim

Taking $Y$ to be a point shows that POPI implies POPA for manifolds.
Taking $Q$ to be empty and $P$ to be a point shows that BOPI implies BOPA.
It also follows that an Oka manifold or a manifold merely satisfying BOPI
has a dominating holomorphic map from Euclidean space; this is not
obvious from the definition of BOPI.

\subhead Basic interpolation implies ellipticity for Stein manifolds
\endsubhead
Ellipticity, introduced by Gromov \cite{G}, and the more general
subellipticity, introduced by Forstneri\v c \cite{F1}, are geometric
sufficient conditions for all Oka properties, including one we have not
mentioned before: the parametric Oka property with {\it jet} interpolation
on a submanifold (POPJI), which involves fixing not just the
values but a finite-order jet of a map along a submanifold as the map is
deformed to a globally holomorphic map.  Gromov noted that for a Stein
manifold, the basic version of this property implies ellipticity
\cite{G, 3.2.A}.  With the help of a mapping cylinder and a more
sophisticated Runge approximation than above we will now show that
ellipticity of a Stein manifold actually follows from the basic Oka
property with plain, zeroth-order interpolation. 

Let $X$ be a Stein manifold.  Without loss of generality we may assume
that $X$ is connected.  The proof of Gromov's result given by
Forstneri\v c and Prezelj \cite{FP, Proposition 1.2} produces a smooth
map $u$ from the tangent bundle $Y=TX$ of $X$, which is a Stein
manifold, to $X$, such that $u$ is the identity on the zero section
$Z$ of $TX$ (when $Z$ is identified with $X$), $u$ is holomorphic on a
neighbourhood $V$ of $Z$, and the derivative of $u$ at a point of $Z$
restricts to the identity map of the tangent space of $X$ at that
point, viewed as a subspace of the tangent space of $Y$ there.  If $u$
was holomorphic on all of $Y$, it would be a dominating spray and $X$
would be elliptic.  The basic Oka property with jet interpolation on
$Z$ allows us to deform $u$ to a dominating spray without further ado;
we will produce a dominating family of sprays on $X$ using only plain
interpolation.

Every analytic subset of a Stein manifold has a fundamental system of
Runge neighbourhoods \cite{KK, E\. 63e}, so we may assume that $V$ is
Runge (as many but not all authors do, we take Runge to include Stein).
Let $\phi:V\to\Bbb C^r$ be a holomorphic
embedding.  The inclusion $\iota:V\hookrightarrow Y$ factors through the
mapping cylinder $M=Y\times\Bbb C^r$ as $V @>j>> M @>\pi>>Y$, where
$j=(\iota,\phi)$ is an embedding and $\pi$ is the projection $(y,z)
\mapsto y$.  Consider the smooth map $u\circ\pi:M\to X$, which is
holomorphic on the submanifold $j(V)$ of $M$.  Assume $X$ satisfies
BOPI.  Then $u\circ\pi$ can be deformed to a holomorphic map $v:M\to X$,
keeping it fixed on $j(V)$.  We plan to precompose $v$ by finitely
many holomorphic sections of $\pi$ of the form $\sigma=(\text{id}_Y,\psi):
Y\to M$, where the holomorphic map $\psi:Y\to\Bbb C^r$ equals $\phi$ on
$Z$ and approximates $\phi$ sufficiently well near a subset of $Z$ in a
sense to be made precise.  Each holomorphic map $v\circ\sigma:Y\to X$ is
the identity on the zero section and its derivative restricted to the
tangent space of $X$ at a point of the zero section is surjective
on the set where $\psi$ approximates $\phi$ well enough.  If these
sets cover the zero section, then these maps form a dominating family
of sprays for $X$, showing that $X$ is subelliptic.  Then, being
Stein, $X$ is in fact elliptic \cite{F1, Lemma 2.2}.

We shall use the method of {\it admissible systems} used by Narasimhan
to prove the embedding theorem for Stein spaces and attributed to Grauert
in \cite{N}.  There are admissible systems $(U^\lambda_i)_{i\in\Bbb N}$,
$\lambda=1,\dots,4\dim X+1$, that together cover $Y$.
The construction in \cite{N} shows that the relatively compact open
subsets $U^\lambda_i$ of $Y$ may be taken to be Stein.  Find compact sets
$K^\lambda_i$ and $L^\lambda_i$ such that $L^\lambda_i\subset
\operatorname{int} K^\lambda_i$, $K^\lambda_i\subset V\cap U^\lambda_i$,
and $Z\subset \bigcup_{\lambda,i}L^\lambda_i$.  Now fix $\lambda$ and
omit it from the notation.  Let $f:V\to\Bbb C$ be one of the coordinate
functions of $\phi$.  We define a holomorphic function $g$ on the union
$Z\cup\bigcup U_i$ as follows.  We set $g=f$ on $Z$.  If $U_i\cap V
=\varnothing$, let $g=0$ on $U_i$.  If $U_i$ intersects $V$, then we take
$g$ on $U_i$ to be an approximation of $f|V\cap U_i$ on $K_i$ with
$g=f$ on $Z\cap U_i$.  Here we apply the Runge Approximation Theorem
for coherent analytic sheaves \cite{Fo, p\. 139} to the ideal sheaf
$\Cal I_Z$ of $Z$ in the Stein manifold $U_i$, using the fact that
$V\cap U_i$ is Runge in $U_i$.  We first extend $f|Z\cap U_i$ to a
holomorphic function $\tilde f$ on $U_i$; then we approximate
$f-\tilde f\in\Cal I_Z(V\cap U_i)$ by $\tilde g\in\Cal I_Z(U_i)$ on
$K_i$ and let $g=\tilde f+\tilde g$.

Now we need the following straightforward adaptation of Narasimhan's
approximation theorem \cite{N, Theorem 2} to the case of an ideal sheaf.

\proclaim{Approximation Theorem}  Let $(U_i)_{i\in\Bbb N}$ be an
admissible system on a Stein manifold $Y$ and $Z$ be an analytic subset
of $Y$.  Let $g$ be a holomorphic function on $Z\cup\bigcup U_i$.  Let
$K_i\subset U_i$ be compact and $\epsilon_i>0$.  Then there is a
holomorphic function $h$ on $Y$ such that $h=g$ on $Z$ and
$|h-g|<\epsilon_i$ on $K_i$ for all $i$.
\endproclaim

We obtain a holomorphic function $h$ on $Y$ with $h=g=f$ on $Z$ such
that $h$ approximates the original function $f$ as closely as
we wish on each compact set $K_i$, so $h$ approximates the
derivatives of $f$ as closely as we wish on each compact set $L_i$.
Approximating each coordinate function of $\phi$ in this way, we obtain
a holomorphic map $\psi:Y\to\Bbb C^r$ such that the composition
$v\circ(\text{id}_Y,\psi):Y\to X$ is the identity on the zero section
and its derivative restricted to the tangent space of $X$ at a point of
$Z\cap\bigcup L_i$ is so close to the derivative of $u$ that it is
surjective.

We have proved the following result.

\proclaim{Theorem 2}  A Stein manifold satisfying the basic Oka property
with interpolation is elliptic.
\endproclaim

Hence, for a Stein manifold, BOPI is equivalent to ellipticity and all
Oka-type properties in between, such as subellipticity, POPI, and POPJI.
In the argument above, we could avoid taking $V$ to be Runge and simplify
the definition of $g$ if we had a detailed proof or a reference to a
proof that a Stein manifold can be covered by finitely many admissible
systems all of which are subordinate to a given open cover.

\subhead Factorization in the projective structure \endsubhead
Turning now from the intermediate structure to the projective structure,
we will explicitly factor an arbitrary holomorphic map $f:R\to S$
between Stein manifolds as a projective cofibration into a
new mapping cylinder $M$ followed by a pointwise weak equivalence.
The new mapping cylinder is constructed as a (pointwise) pushout 
$$\xymatrix{
R\ar[r]^f \ar[d]^{i_0} &  S \ar[d] \\
R\times\Delta^1 \ar[r] & M
}$$
in $\frak S$, or equivalently as the pushout
$$\xymatrix{
R\sqcup R \ar[r]^{f\sqcup \text{id}} \ar[d]^{i_0\sqcup i_1} &
S\sqcup R \ar[d] \\ R\times\Delta^1 \ar[r] & M
}$$
Since the left-hand vertical maps are cofibrations in the projective
structure, so are the right-hand vertical maps.  Since $S$ is
cofibrant, the map $R\to S\sqcup R$ that is the identity on $R$ is a
cofibration, so the composition $j:R\to S\sqcup R \to M$ is a cofibration.
Since $R$ is cofibrant, it follows that $M$ is cofibrant as well.
The original map $f:R\to S$ factors as $\pi\circ j$,
where $\pi$ is the obvious projection $M\to S$, which is a left inverse to
the map $S\to M$ given directly by the construction of $M$.  Now $R\to
R\times\Delta^1$ is a pointwise acyclic cofibration, so its pushout
$S\to M$ is pointwise acyclic.  Hence, the left inverse $\pi$ is pointwise
acyclic as well.  It follows that $j$ is acyclic if and only if $f$ is
acyclic, that is, $f$ is a topological homotopy equivalence.

A map in $\frak S$ from $M$ to a manifold $X$ is easy to understand
explicitly.  It is nothing but a pair of maps, one from
$\hom_{\frak S}(S,X)=\Cal O(S,X)$ and another from
$$\hom_{\frak S}(R\times\Delta^1,X)=\hom_{s\text{\bf
Set}}(\Delta^1,\Hom_{\frak S}(R,X))=s\Cal O(R,X)_1=\Cal C(I,\Cal
O(R,X))$$
(here, $s\text{\bf Set}$ denotes the category of simplicial sets),
which agree when the first map is precomposed with
$f:R\to S$ and the second is restricted to $\{0\}\subset I$ (the pushout
that defines $M$ is taken levelwise as well as pointwise).  In other
words, a map $M\to X$ is nothing but a continuous map from the ordinary
topological mapping cylinder $S\cup_f(R\times I)$ to $X$ which
is holomorphic on each slice $S$ and $R\times\{t\}$, $t\in I$.  To put
it a third way, a map $M\to X$ is simply a continuous deformation of a
holomorphic map $R \to X$ through such maps to a map that factors
through $f$.

\subhead An application \endsubhead
If $f$ is a topological weak equivalence, so $j$ is an acyclic
cofibration, and $X$ is weakly Oka, it follows that the map
$\Hom(M,X)\to\Hom(R,X)=s\Cal O(R,X)$ induced by $j$ is an acyclic
fibration and hence surjective on vertices.  This means that every
holomorphic map $R\to X$ can be continuously deformed through such maps
to a map that factors through $f$ by a holomorphic map $S\to X$.
(This could for instance be applied to zero-free holomorphic
functions, taking $X=\Bbb C\setminus\{0\}$, and to meromorphic functions
without indeterminacies, taking $X$ to be the Riemann sphere.)  The same
holds if $f$ is arbitrary and $X$ is weakly Oka and topologically
contractible.

\subhead The topological mapping cylinder \endsubhead
Let us finally apply the ordinary topological mapping cylinder $M$ of
the inclusion of an open Stein subset $R$ in a Stein manifold $S$,
defined as the subset of $S\times[0,1]$ of pairs $(x,t)$ for which $x\in
R$ if $t>0$.  The inclusion $R\hookrightarrow S$ factors through $M$ as the
inclusion $R\hookrightarrow M$, $x\mapsto (x,1)$, which is a topological
cofibration (in the stronger of the two senses used by topologists,
namely the sense that goes with Serre fibrations rather than Hurewicz
fibrations), followed by the projection $M\to S$, $(x,t)\mapsto x$, which
is a topological homotopy equivalence.  

Let $X$ be a weakly Oka manifold.  Recall that this means that the inclusion
$\Cal O(T,X) \hookrightarrow\Cal C(T,X)$ is a weak equivalence in the
compact-open topology for every Stein manifold $T$.  Then we have the
following diagram.
$$\xymatrix{
\Cal O(S,X) \ar[dd]^{\text{w.eq.}} \ar[rrd] \ar@{.>}[rd]_{\text{w.eq.}} & & \\
& \text{pullback} \ar[d]^{\text{w.eq.}} \ar[r] & \Cal O(R,X) \ar[d]^{\text{w.eq}} \\
\Cal C(S,X) \ar [r]^{\text{w.eq.}} & \Cal C(M,X) \ar[r]^{\text{fibr.}} & \Cal C(R,X)
}$$
The left-hand and right-hand vertical maps are weak equivalences because
$X$ is weakly Oka.  The right-hand bottom map is a fibration because
$R\hookrightarrow M$ is a cofibration.  The middle
vertical map is a weak equivalence because it is the pullback of a weak
equivalence by a fibration.  The left-hand bottom map is a weak
equivalence because $M\to S$ is a homotopy equivalence.  Finally, we conclude
that the dotted map
$$\Cal O(S,X) \to \Cal C(M,X)\times_{\Cal C(R,X)}\Cal O(R,X)$$
is a weak equivalence.

\subhead One more implication between Oka properties \endsubhead
The dotted map being a weak equivalence implies that a continuous map $S\to X$
which is holomorphic on $R$ can be deformed through continuous maps $M\to X$
which are holomorphic on $R\subset M$ to a holomorphic map $S\to X$.  Take
open subsets $U$ and $V$ of $R$ with $\overline U\subset V$ and $\overline
V\subset R$.  Precomposing the deformation by a section of the projection
$M\to S$ of the form $x\mapsto (x,\rho(x))$, where $\rho:S\to[0,1]$ is
continuous, $\rho=1$ on $U$, and $\rho=0$ outside $V$, shows that a
continuous map $S\to X$ which is holomorphic on $R$ can be deformed through
continuous maps $S\to X$ which are holomorphic on $U$ to a holomorphic map
$S\to X$.

Thus wanting to keep a continuous map $S\to X$ holomorphic on a Stein
open set where it happens to be holomorphic to begin with, as it is deformed
to a holomorphic map, does not lead to a new version of the Oka property,
at least as long as we do not mind slightly shrinking the open subset.  (This
is in the same spirit as our earlier result that POPI subsumes POPA.)  It is
not known if shrinking is actually necessary.

\subhead Remarks \endsubhead
Our results so far clarify and simplify the spectrum of Oka
properties and their significance.  It is beginning to seem plausible
that there is essentially only one Oka property, the one called {\it the}
Oka property in \cite{L2}, and that other Oka properties should be viewed
as variants or special cases of it.  (Although they are technically more
difficult to work with, the parametric properties are more
interesting than the basic ones from a homotopy-theoretic point of view:
it is more natural to work with weak equivalences of mapping spaces
than maps that merely induce surjections of path components.)  It also
appears that there is no reason to look for intermediate model structures
of relevance to Oka theory beyond the one already defined in \cite{L2}.

Furthermore, the chief practical sufficient condition for the Oka
property, ellipticity, turns out to be equivalent to it, at least for
Stein manifolds.  To what extent this generalizes to arbitrary manifolds
is a most interesting open question.  If $X$ is any manifold satisfying
BOPI (or just BOPA) and $p\in X$, then there is a holomorphic
map from Euclidean space of dimension $\dim_p X$ to $X$ that takes the 
origin to $p$ and is a submersion there.  The question is whether such 
maps can be made to fit together to make a dominating spray or a 
dominating family of sprays or some weaker geometric structure that
could still be used to establish the Oka property.  In the remainder of
the paper we take a step towards answering this question in the
affirmative.

\subhead Cofibrant models \endsubhead
A cofibrant model for a complex manifold $X$ in the
intermediate model structure on $\frak S$ is an acyclic
intermediate fibration from an intermediately cofibrant prestack
to $X$, or, in other words, a factorization of the map from
the empty manifold to $X$ of the same sort as before.  (There
is a distinction between the empty prestack (the initial object
in $\frak S$) and the prestack represented by the empty manifold,
but here it is immaterial.)  Such a factorization always
exists by the axioms for a model category, but we seek a
reasonably explicit cofibrant model represented by a manifold $S$.  An
intermediately cofibrant manifold is Stein --- we defer the rather
technical proof of this to the end of the paper --- so $S$ will be
Stein and the acyclic intermediate fibration $S\to X$ holomorphic.

Our motivation for seeking cofibrant models represented by 
Stein manifolds is threefold.  First, since fibrancy passes up and down
in an acyclic fibration, a Stein cofibrant model allows us to reduce
the problem of geometrically characterizing the Oka property to the
Stein case, of which we have already disposed.  Second, an acyclic
intermediate fibration from a Stein manifold to $X$ is weakly final
among maps from Stein manifolds into $X$: all such maps factor through
it (not necessarily uniquely).  This is a nontrivial and apparently new
notion that may be of independent complex-analytic interest.  Third,
from a homotopy-theoretic point of view, we are taking explicit
factorizations beyond the relatively easy case in which the source and
target are both cofibrant.

\subhead The Jouanolou Trick \endsubhead
The so-called Jouanolou Trick, invented for the
purpose of extending $K$-theory from affine schemes to quasi-projective
ones \cite{J}, is the observation that every quasi-projective
scheme $X$ carries an affine bundle whose total space $A$ is affine.
Then $X$ and $A$ have the same motive.  This can be used to reduce
various questions in algebraic geometry to the affine case.  Of interest
here is the expectation that $A$ is a cofibrant model for $X$
in the intermediate structure.  Let us outline an analytic version of
the Jouanolou Trick, starting with projective space.

Let $\Bbb P_n$ denote $n$-dimensional complex projective space.  Let
$\Bbb Q_n$ be the complement in $\Bbb P_n\times\Bbb P_n$ of the
hypersurface of points $([z_0:\dots:z_n],[w_0:\dots:w_n])$ with
$z_0 w_0+\dots+z_n w_n=0$.  This hypersurface is the preimage of a
hyperplane by the Segr\acuteaccent e embedding $\Bbb P_n\times
\Bbb P_n\to \Bbb P_{n^2+2n}$, so $\Bbb Q_n$ is Stein.  Let $\pi$ be
the projection $\Bbb Q_n\to\Bbb P_n$ onto the first component.
It is easily seen that $\pi$ has the structure of an
affine bundle with fibre $\Bbb C^n$ (of course without holomorphic
sections).   In particular, $\pi$ is an acyclic elliptic
bundle and hence an acyclic subelliptic submersive Serre fibration
(SSSF).  Our conjecture from \cite{L2} that an SSSF is an intermediate
fibration remains open, so for now $\pi$ is only a candidate for a
cofibrant model of $\Bbb P_n$.  We shall continue with SSSFs standing
in for intermediate fibrations.

\subhead Towards the general case \endsubhead
Choose a class $\Cal G$ of acyclic SSSFs, containing all affine bundles
and closed under taking pullbacks by arbitrary holomorphic maps.
Call the maps in $\Cal G$ {\it good} maps.  For example, $\Cal G$ could
be the class of affine bundles; holomorphic fibre bundles with a
contractible subelliptic fibre, such as $\Bbb C^k$ for some $k\geq 0$;
locally smoothly or real-analytically trivial acyclic subelliptic
submersions; or all acyclic SSSFs.  Call a complex manifold {\it good}
if it is the target (and hence the image) of a good map from a Stein
manifold.

We have seen that projective spaces are good and so are Stein manifolds,
obviously.  The pullback of a good map $g:S\to Y$ by any holomorphic
map $f:X\to Y$ is a good map $f^*g:R\to X$ ($R$ is smooth, that is,
a manifold, because $g$ is a submersion).  If $f$ is a covering map, a
finite branched covering map, the inclusion of a (closed) submanifold,
or the inclusion of the complement of an analytic hypersurface, then so
is the pullback map $g^*f:R\to S$, so if $S$ is Stein, then $R$ is also
Stein.  By blowing up the complement of a Zariski-open subset of a
projective variety, turning it into a hypersurface, we see that
quasi-projective manifolds are good.  Finally, it is easy to see that
the product of good manifolds is good.

The class of good manifold thus appears to be quite large (even with
the smallest possible $\Cal G$).  It contains all Stein manifolds and
quasi-projective manifolds and is closed under taking products,
covering spaces, finite branched covering spaces,
submanifolds, and complements of analytic hypersurfaces.  We do not
know if every manifold, or even every domain in Euclidean space, is the
target of an acyclic SSSF from a Stein manifold.

An acyclic SSSF from a Stein manifold $S$ to $X$ is weakly final among
holomorphic maps from Stein manifolds to $X$.  Namely, if $R$ is Stein
and $R\to X$ is a holomorphic map, then there is a continuous lifting
$h:R\to S$ since $S\to X$ is an acyclic Serre fibration, and because
$S\to X$ is a subelliptic submersion, $h$ can be deformed to a
holomorphic lifting.  (This can of course be strengthened to include
interpolation on a submanifold of $R$.)  We do not know if every
manifold has a weakly final map from a Stein manifold.

\smallskip
The next theorem is a weak variant of one of the interesting statements
that follow from the conjecture that an SSSF is an intermediate
fibration.

\proclaim{Theorem 3}  Let $X$ and $Y$ be complex manifolds and $f:X\to
Y$ be a submersive subelliptic Serre fibration.  If $Y$ has the basic
Oka property with interpolation, then so does $X$.  If $f$ is acyclic
and $X$ has the basic Oka property with interpolation, then so does $Y$.
\endproclaim

This statement with BOPI replaced by POPI is a direct consequence of
our conjecture since POPI is equivalent to intermediate fibrancy
\cite{L2}, but BOPI is all we can handle at present.  Forstneri\v c
has proved an analogous result for BOPA \cite{F3}.  (He calls BOPA
{\it the Oka property}.)  Interestingly enough, he does not require
acyclicity in order to pass from the source to the target.  We do not
know if this holds for BOPI.  It is also proved in \cite{F3} that POPA
passes up in an SSSF.

\demo{Proof}  First suppose $Y$ satisfies BOPI.  Let $h:S\to X$ be a
continuous map from a Stein manifold $S$ such that the restriction of $h$
to a submanifold $T$ of $S$ is holomorphic.  We need to show that $h$
can be deformed to a holomorphic map $S\to X$ keeping it fixed on $T$.

By assumption, the composition $f\circ h$ can be deformed to a
holomorphic map $g:S\to Y$ keeping it fixed on $T$.  By the topological
Axiom SM7, since $f$ is a Serre fibration, the map
$$\Cal C(S,X)\to\Cal C(S,Y)\times_{\Cal C(T,Y)}\Cal C(T,X)$$
is a fibration, so as $f\circ h$ is deformed to $g$, we can deform $h$
with it to a continuous map $k:S\to X$, keeping it fixed on $T$.  This
yields the following square of holomorphic maps.
$$\xymatrix{
T \ar[r]\ar[d] & X \ar[d]^{f} \\
S \ar[r]^{g} \ar@{-->}[ur]^{k} & Y
}$$
Since $f$ is a subelliptic submersion, we can deform $k$ to a
holomorphic lifting in this square.  Thus, in two steps, $h$ has been
deformed to a holomorphic map, keeping it fixed on $T$.

Second, suppose $f$ is acyclic and $X$ satisfies BOPI.  Let $h:S\to Y$
be a continuous map from a Stein manifold $S$ such that the restriction
of $h$ to a submanifold $T$ of $S$ is holomorphic.  Since $f$ is an acyclic
Serre fibration, there is a continuous lifting $T\to X$ of $h|T$.
Since $f$ is a subelliptic submersion, this lifting can be deformed to a
holomorphic lifting.  Again, since $f$ is an acyclic Serre fibration,
there is a continuous lifting in the square
$$\xymatrix{
T \ar[r]^{\text{hol\.}} \ar[d] & X \ar[d]^{f} \\
S \ar[r]^{h} \ar@{-->}[ur] & Y
}$$
Since $X$ satisfies BOPI, the lifting can be deformed to a holomorphic
map, keeping it fixed on $T$.  Postcomposing with $f$ gives a
deformation of $h$ to a holomorphic map, keeping it fixed on $T$, and
the proof is complete.
\qed\enddemo

Taking $T$ to be empty proves the theorem with BOPI replaced by BOP.
For BOP, acyclicity is required in order to pass from the source to the
target, as shown by the simple example of the universal covering map
from the disc to the punctured disc.

A holomorphic intermediate fibration satisfies the conclusion of
Theorem 3, because, as shown in \cite{L2}, it has the properties
of an SSSF that were used in the proof above.

It is an open question whether a manifold $X$ fibred over an elliptic
manifold by elliptic manifolds is elliptic.  It follows from Theorems 2
and 3 that the answer is affirmative if $X$ is Stein; in general $X$ at
least satisfies BOPI by Theorem 3.

\smallskip
Here is our geometric characterization of the basic Oka
property with interpolation for an arbitrary good manifold, extending
the characterization by ellipticity for Stein manifolds given by Theorem
2.

\proclaim{Theorem 4}  A good manifold has the basic Oka property with
interpolation if and only if it is the image of a good map from an
elliptic Stein manifold.
\endproclaim

\demo{Proof}  Let $X$ be a manifold and $S\to X$ be a good map.
If $S$ is elliptic, then $S$ satisfies BOPI, so $X$ does too by
Theorem 3.  (For this direction, we do not need $S$ to be Stein.)
Conversely, suppose $X$ satisfies BOPI and is good, so there is a
good map $S\to X$ with $S$ Stein.  Then $S$ satisfies BOPI by
Theorem 3 and is therefore elliptic by Theorem 2. 
\qed\enddemo

We can state this more explicitly for quasi-projective manifolds.

\proclaim{Theorem 5}  A quasi-projective manifold has the basic Oka
property with interpolation if and only if it carries an affine bundle
whose total space is elliptic and Stein.
\endproclaim

Note that if good maps are intermediate fibrations, as
predicted by our conjecture, so POPI passes down in a good map,
then it follows that BOPI and POPI are equivalent for good manifolds.

\subhead A converse Oka Principle \endsubhead
By the definition of the intermediate model structure, a Stein
manifold is intermediately cofibrant.  We conclude this paper by proving
that, conversely, an intermediately cofibrant complex manifold is Stein.
It is then immediate that projective cofibrancy is also equivalent to
being Stein.  This may be viewed as a \lq\lq converse Oka
Principle\rq\rq\ for manifolds.  Characterizing the Stein property as
cofibrancy complements the characterizations in \cite{L2} of the
parametric Oka property with and without interpolation as intermediate
and projective fibrancy, respectively. We hope to treat the general
case of holomorphic maps later.

\subhead Remark on fullness \endsubhead
Before proceeding to the proof, we recall from \cite{L2} that while the
Yoneda embedding of the category $\Cal S$ of Stein manifolds into the
category $\frak S$ of prestacks on $\Cal S$ is of course full (even
simplicially full), we have been unable to determine whether its
extension to an embedding of the category $\Cal M$ of all complex
manifolds into $\frak S$, taking a manifold $Z$ to the prestack
$s\Cal O(\cdot, Z)$, is full.  It is easy to demonstrate a
weaker fullness property, which is required here and suffices for many
purposes.  A morphism $\eta:X\to Y$ in $\frak S$ between complex
manifolds (that is, between the prestacks represented by them)
is a natural transformation $s\Cal O(\cdot,X)\to s\Cal O(\cdot,Y)$ of
contravariant simplicial functors $\Cal S \to s\text{\bf Set}$.
Considering what this means, see e.g\. the definitions in
\cite{GJ, Chapter IX}, we observe that $\eta$ {\it restricted to the
vertex level} is nothing but a morphism $\Cal O(\cdot, X)\to
\Cal O(\cdot, Y)$ in the category $\operatorname{Pre}\Cal S$ of
presheaves of sets on $\Cal S$.  Such a morphism is induced by a
holomorphic map $X\to Y$ since by \cite{L1, Proposition 4.2}, the
functor $\Cal M\to \operatorname{Pre}\Cal S$ taking a manifold $Z$ to
$\Cal O(\cdot, Z)$ is a full embedding.  What we do not know unless $X$
is Stein is whether this holomorphic map also induces the higher-level
maps in $\eta$.  For all we know, distinct morphisms $X\to Y$ in
$\frak S$ might coincide at the vertex level.

We conclude that if we restrict any diagram of complex manifolds and
morphisms in $\frak S$ to the vertex level, the morphisms will be
induced by holomorphic maps, so we get a diagram in $\Cal M$.  In other
words, restriction to the vertex level gives a retraction functor from
the full subcategory of $\frak S$ generated by $\Cal M$ onto $\Cal M$
itself.

\proclaim{Theorem 6}  Let $X$ be a complex manifold.  The following are
equivalent.
\roster
\item $X$ is intermediately cofibrant.
\item $X$ is projectively cofibrant.
\item $X$ is Stein.
\endroster
\endproclaim

\demo{Proof}  By the definition of the intermediate and projective model
structures on $\frak S$, it is clear that \therosteritem3 $\Rightarrow$
\therosteritem2 $\Rightarrow$ \therosteritem1.

Assuming that $X$ is intermediately cofibrant, we need to prove that
$X$ is Stein.  It suffices to show that we can prescribe the values
of a holomorphic function on $X$ on any two-point subset of $X$ (this
gives holomorphic separability) and on some infinite subset of any
infinite discrete subset of $X$ (this gives holomorphic convexity).  By
definition of the intermediate model structure, $X$ is a
retract $X\to Y\to X$ of a cell complex $Y$, that is, a prestack which
is the target of a transfinite composition with source $\varnothing$ of
pushouts of generating cofibrations
$$S\times\partial\Delta^n \cup_{T\times\partial\Delta^n} T\times\Delta^n
\to S\times\Delta^n,$$
where $T\hookrightarrow S$ is a Stein inclusion and $n\geq 0$.  More
explicitly, $Y=\colim_{\alpha<\lambda} A_\alpha$, where $A:\lambda\to
\frak S$ is a functor from an ordinal $\lambda$ with $A_0=\varnothing$
such that for every limit ordinal $\gamma<\lambda$, the induced map
$\colim_{\alpha<\gamma} A_\alpha\to A_\gamma$ is an isomorphism, and
such that for every successor ordinal $\alpha<\lambda$, the map
$A_{\alpha-1} \to A_\alpha$ is a pushout of a generating cofibration
given by a Stein inclusion $T_\alpha\hookrightarrow S_\alpha$ and an
integer $n_\alpha\geq 0$.

Let $E$ be an infinite discrete subset of $X$, viewed as a discrete
Stein manifold.  Consider the subset of $\lambda$ of ordinals $\alpha$
such that the restriction of $E\hookrightarrow X\to Y$ to an infinite
subset of $E$ factors through $S_\alpha\times\Delta^{n_\alpha}$.  This
subset is not empty: since colimits in $\frak S$ are taken levelwise
and pointwise, the map $E\to Y$ itself, as a vertex of $Y(E)$, is the
image of a vertex of $(S_\alpha\times\Delta^{n_\alpha})(E)$ for some
$\alpha$.  Let $\beta$ be the smallest element of this subset and let
$E'$ be an infinite subset of $E$ such that the restriction $E'\to Y$
factors through $S_\beta\times\Delta^{n_\beta}$.  At the vertex level,
generating cofibrations with $n\geq 1$ are isomorphisms since
$\partial\Delta_0^n = \{0,\dots,n\} = \Delta_0^n$, so $n_\beta=0$, for
otherwise $E'\to Y$ would factor through $A_{\beta-1}$ and hence
through $S_\alpha\times\Delta^{n_\alpha}$ for some $\alpha<\beta$.
Since $E'\to Y$ is a monomorphism, so is $E'\to S_\beta$.  The
composition $E'\to S_\beta\to Y\to X$, which is simply the inclusion
of $E'$ in $X$, has discrete image, and the morphism $S_\beta\to X$ is
induced by a holomorphic map, so $E'\to S_\beta$ has discrete image
and we can view $E'$ as a discrete subset of $S_\beta$.  The
restriction of $E\to Y$ to $E'\cap T_\beta$ factors through
$A_{\beta-1}$ and hence through $S_\alpha\times\Delta^{n_\alpha}$ for
some $\alpha<\beta$, so $E'\cap T_\beta$ is finite.  Let $F$ be the
infinite subset $E'\setminus T_\beta$ of $E$.

Similarly, if $F$ is a two-point subset of $X$, let $\beta$ be the
smallest element of $\lambda$ such that $F\hookrightarrow X\to Y$
factors through $S_\beta\times\Delta^{n_\beta}$.  Then $n_\beta=0$
and $T_\beta$ cannot contain $F$.

In either of the two cases, take a map $F\to\Bbb C$ and extend it to
$F\cup T_\beta$ by mapping $T_\beta$ to a constant.
Map $A_{\beta-1}$ to the same constant, say through the map to the
final prestack.  Since $F\cup T_\beta$ is a submanifold of the
Stein manifold $S_\beta$, the holomorphic map $F\cup T_\beta\to\Bbb C$ 
has a holomorphic extension to $S_\beta$.
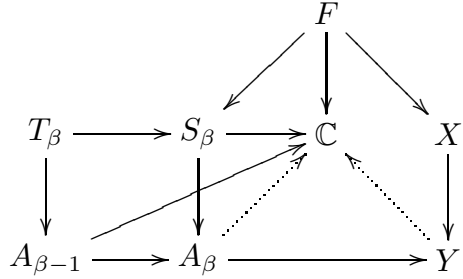
$$\xymatrix{
 & & F \ar[dl] \ar[d] \ar[dr] & \\
T_\beta \ar[r] \ar[d] & S_\beta \ar[d] \ar[r] & \Bbb C & X \ar[d] \\
A_{\beta-1} \ar[urr] \ar[r] & A_\beta \ar[rr] \ar@{.>}[ur] & &
Y \ar@{.>}[ul]
}$$
By the pushout property of $A_\beta$, we get an induced map
$A_\beta\to\Bbb C$.  Since $A_\beta\to Y$ is an intermediate
cofibration and $\Bbb C$ is intermediately fibrant and acyclic,
there is an extension $Y\to\Bbb C$.  Finally, the restriction to the
vertex level of the composition $X\to Y\to \Bbb C$ is a holomorphic
map by our fullness remark and it extends the given map on $F$.
\qed\enddemo

\enddocument